\documentclass[reqno,centertags]{amsart}
\usepackage{eurosym}
\usepackage{amsmath}
\usepackage{amssymb}
\usepackage{amsfonts}

\newtheorem{theorem}{Theorem}
\theoremstyle{plain}

\newtheorem{corollary}{Corollary}

\newtheorem{lemma}{Lemma}

\newtheorem{problem}{Problem}

\newtheorem{remark}{Remark}

\numberwithin{equation}{section}

\input{tcilatex}

\begin{document}
\title[]{Tubes of finite Chen-type}
\author{Hassan Al-Zoubi }
\address{Department of Mathematics, Al-Zaytoonah University of Jordan, P.O.
Box 130, Amman, Jordan 11733}
\email{dr.hassanz@zuj.edu.jo}
\author{Stylianos Stamatakis}
\address{Department of Mathematics, Aristotle University of Thessaloniki}
\email{stamata@math.auth.gr}
\author{Khalid M. Jaber}
\address{Al-Zaytoonah University of Jordan, Department of Computer Science }
\email{K.Jaber@zuj.edu.jo}
\author{Hani Almimi}
\address{Al-Zaytoonah University of Jordan, Department of Computer Science}
\email{Hani.Mimi@zuj.edu.jo}
\date{}
\subjclass[2010]{ 53A05}
\keywords{ Surfaces in the Euclidean 3-space, Surfaces of finite Chen-type,
Beltrami operator. }

\begin{abstract}
In this paper, we consider surfaces in the 3-dimensional Euclidean space $%
E^{3}$ which are of finite $III$-type, that is, they are of finite type, in
the sense of B.-Y. Chen, corresponding to the third fundamental form. We
present an important family of surfaces, namely, tubes in $E^{3}$. We show
that tubes are of infinite $III$-type.
\end{abstract}

\maketitle

\section{Introduction}

Let $M^{n}$ be a (connected) submanifold in the m-dimensional Euclidean
space $E^{m}$. Let $\boldsymbol{x},\boldsymbol{H}$ be the position vector
field and the mean curvature field of $M^{n}$ respectively. Denote by $%
\Delta ^{I}$ the second Beltrami-Laplace operator corresponding to the first
fundamental form $I$ of $M^{n}$. Then, it is well known that \cite{R3}
\begin{equation}
\Delta ^{I}\boldsymbol{x}=-n\boldsymbol{H}.  \notag
\end{equation}

From this formula one can see that $M^{n}$ is a minimal submanifold if and
only if all coordinate functions, restricted to $M^{n}$, are eigenfunctions
of $\Delta ^{I}$ with eigenvalue $\lambda =0$. Moreover in \cite{R15} T.
Takahashi showed that the submanifold $M^{n}$ for which $\Delta ^{I}%
\boldsymbol{x}=\lambda \boldsymbol{x}$, i.e., for which all coordinate
functions are eigenfunctions of $\Delta ^{I}$ with the same eigenvalue $%
\lambda \in \mathbb{Re}$, are precisely either the minimal submanifold with
eigenvalue $\lambda =0$ or the minimal submanifold of hyperspheres $S^{m-1}$
with eigenvalue $\lambda >0$. Although the class of finite type submanifolds
in an arbitrary dimensional Euclidean spaces is very large, very little is
known about surfaces of finite type in the Euclidean 3-space $E^{3}$.
Actually, so far, the only known surfaces of finite type corresponding to
the first fundamental form in the Euclidean 3-space are the minimal
surfaces, the circular cylinders and the spheres. So in \cite{R5} B.-Y. Chen
mentions the following problem

\begin{problem}
\label{(1)}Determine all surfaces of finite Chen $I$-type in $E^{3}$.
\end{problem}

In order to provide an answer to the above problem, important families of
surfaces were studied by different authors by proving that finite type ruled
surfaces, finite type quadrics, finite type tubes, finite type cyclides of
Dupin and finite type spiral surfaces are surfaces of the only known
examples in $E^{3}$. However, for another classical families of surfaces,
such as surfaces of revolution, translation surfaces as well as helicoidal
surfaces, the classification of its finite type surfaces is not known yet.
For a more details, the reader can refer to \cite{R6}.

Later in \cite{R11} O. Garay generalized T. Takahashi's condition studied
surfaces in $E^{3}$ for which all coordinate functions $\left(
x_{1},x_{2},x_{3}\right)$ of $\boldsymbol{x}$ satisfy $\Delta^{I}\boldsymbol{%
x}_{i} = \lambda_{i}x_{i}, i = 1,2,3$, not necessarily with the same
eigenvalue. Another generalization was studied in \cite{R9} for which
surfaces in $E^{3}$ satisfy the condition $\Delta^{I}\boldsymbol{x}= A%
\boldsymbol{x} + B$ $(\ddag)$ where $A \in\mathbb{Re}^{3\times3} ;B \in
\mathbb{Re}^{3\times1}$. It was shown that a surface $S$ in $E^{3}$
satisfies $(\ddag)$ if and only if it is an open part of a minimal surface,
a sphere, or a circular cylinder. Surfaces satisfying $(\ddag)$ are said to
be of coordinate finite type.

In the thematic circle of the surfaces of finite type in the Euclidean space
$E^{3}$, S. Stamatakis and H. Al-Zoubi in \cite{R13} restored attention to
this theme by introducing the notion of surfaces of finite type
corresponding to the second or the third fundamental forms of $S$ in the
following way:

A surface $S$ is said to be of finite type corresponding to the fundamental
form $J$, or briefly of finite $J$-type, $J=II,III$, if the position vector $%
\boldsymbol{x}$ of $S$ can be written as a finite sum of nonconstant
eigenvectors of the operator $\Delta ^{J}$, that is if
\begin{equation}
\boldsymbol{x}=\boldsymbol{x}_{0}+\sum_{i=1}^{k}\boldsymbol{x}_{i},\ \ \ \ \
\ \Delta ^{J}\boldsymbol{x}_{i}=\lambda _{i}\boldsymbol{x}_{i},\ \ \
i=1,...,k,  \label{0}
\end{equation}%
where $\boldsymbol{x}_{0}$ is a fixed vector and $\boldsymbol{x}_{1},...,%
\boldsymbol{x}_{k}$ are nonconstant maps such that $\Delta ^{J}\boldsymbol{\
x}_{i}=\lambda _{i}\boldsymbol{x}_{i},i=1,...,k$. If, in particular, all
eigenvalues $\lambda _{1},\lambda _{2},...,\lambda _{k}$ are mutually
distinct, then $S$ is said to be of $J$-type $k$, otherwise $S$ is said to
be of infinite type. When $\lambda _{i}=0$ for some \emph{i} = 1,..., \emph{k%
}, then $S$ is said to be of null $J$-type $k$.

In general when $S$ is of finite type $k$, it follows from (\ref{0}) that
there exist a monic polynomial, say $R(x)\neq 0,$ such that $R(\Delta ^{J})(%
\boldsymbol{x}-\boldsymbol{c})=\mathbf{0}.$ Suppose that $R(x)=x^{k}+\sigma
_{1}x^{k-1}+...+\sigma _{k-1}x+\sigma _{k},$ then coefficients $\sigma _{i}$%
\ are given by

\begin{eqnarray}
\sigma _{1} &=&-(\lambda _{1}+\lambda _{2}+...+\lambda _{k}),  \notag \\
\sigma _{2} &=&(\lambda _{1}\lambda _{2}+\lambda _{1}\lambda_{3}+...+\lambda
_{1}\lambda _{k}+\lambda _{2}\lambda _{3}+...+\lambda _{2}\lambda
_{k}+...+\lambda _{k-1}\lambda _{k}),  \notag \\
\sigma _{3} &=&-(\lambda _{1}\lambda _{2}\lambda _{3}+...+\lambda
_{k-2}\lambda _{k-1}\lambda _{k}),  \notag \\
&&.............................................  \notag \\
\sigma _{k} &=&(-1)^{k}\lambda _{1}\lambda _{2}...\lambda _{k}.  \notag
\end{eqnarray}

Therefore the position vector $\boldsymbol{x}$ satisfies the following
equation, (see \cite{R3})

\begin{equation}  \label{1}
(\Delta^{J})^{k}\boldsymbol{x}+\sigma_{1}(\Delta^{J})^{k-1}\boldsymbol{x}%
+...+\sigma_{k}( \boldsymbol{x}-\boldsymbol{c})=\boldsymbol{0}.
\end{equation}

In this paper we will pay attention to surfaces of finite $III$-type.
Firstly, we will establish a formula for $\Delta^{III}\boldsymbol{x}$ by
using Cartan's method of the moving frame. Further, we continue our study by
proving finite type surfaces for an important class of surfaces, namely,
tubes in $E^{3}$.

\section{Preliminaries}

Let $S$ be a (connected) surface in the Euclidean 3-space $E^{3}$, whose
Gaussian curvature $K$ never vanishes. Let $\wp =\{\boldsymbol{\varepsilon
_{1}}(u,v),\boldsymbol{\varepsilon _{2}}(u,v),\boldsymbol{\varepsilon _{3}}%
(u,v)\}$ is a moving frame of $S$, $\boldsymbol{\varepsilon _{3}}=%
\boldsymbol{n}$ is the Gauss map of $S$ and $\det (\boldsymbol{\varepsilon
_{1}},\boldsymbol{\varepsilon _{2}},\boldsymbol{\varepsilon _{3}})=1$. Then
it is well known that there are linear differential forms $\omega
_{1},\omega _{2},\omega _{31},\omega _{32}$ and $\omega _{12}$, such that
\cite{R10}

\begin{equation*}
d\boldsymbol{x}=\omega _{1}\boldsymbol{\varepsilon _{1}}+\omega _{2}%
\boldsymbol{\varepsilon _{2}},\ \ \ d\boldsymbol{n}=\omega _{31}\boldsymbol{%
\varepsilon _{1}}+\omega _{32}\boldsymbol{\varepsilon _{2}},
\end{equation*}

\begin{center}
\begin{equation*}
d\boldsymbol{\varepsilon _{1}}=\omega _{12}\boldsymbol{\varepsilon _{2}}%
-\omega _{31}\boldsymbol{\varepsilon _{3}},\ \ \ d\boldsymbol{\varepsilon
_{2}}=-\omega _{12}\boldsymbol{\varepsilon _{1}}-\omega _{32}\boldsymbol{%
\varepsilon _{3}},
\end{equation*}
\end{center}

and functions $a,b,c,q_{1},q_{2}$ of $u,v$ such that

\begin{equation*}
\omega _{31}=-a\omega _{1}-b\omega _{2},\ \ \ \omega _{32}=-b\omega
_{1}-c\omega _{2},\ \ \ \omega _{12}=q_{1}\omega _{1}+q_{2}\omega _{2}.
\end{equation*}

We can choose the moving frame of $S$, such that the vectors $\boldsymbol{%
\varepsilon_{1}},\boldsymbol{\varepsilon_{2}}$ are the principle directions
of $S$. Then $a$, $c$ are the principle curvatures of $S$ and $b=0$, so the
differential forms $\omega_{1}$ and $\omega_{2}$ become

\begin{equation}
\omega_{1} =-\frac{1}{a}\omega_{31}, \ \ \ \ \omega_{2} =-\frac{1}{c}%
\omega_{32}.  \notag
\end{equation}
The Gauss and mean curvature are respectively

\begin{equation}
K= ac, \ \ \ \ H=\frac{a+c}{2}.  \notag
\end{equation}

Let $\nabla _{1}f,\nabla _{2}f$ be the derivatives of Pfaff of a function $%
f(u,v)\in C^{1}$ along the curves $\omega _{2}=0,\omega _{1}=0$
respectively. Then we have the following well known relations \cite{R2}

\begin{eqnarray*}
\nabla _{1}\boldsymbol{x} &=&\boldsymbol{\varepsilon _{1}},\ \ \ \nabla _{2}%
\boldsymbol{x}=\boldsymbol{\varepsilon _{2}}, \\
\nabla _{1}\boldsymbol{\varepsilon _{1}} &=&q_{1}\boldsymbol{\varepsilon _{2}%
}+a\boldsymbol{n},\ \ \ \nabla _{2}\boldsymbol{\varepsilon _{1}}=q_{2}%
\boldsymbol{\varepsilon _{2}}, \\
\nabla _{1}\boldsymbol{\varepsilon _{2}} &=&-q_{1}\boldsymbol{\varepsilon
_{1}},\ \ \ \ \nabla _{2}\boldsymbol{\varepsilon _{2}}=-q_{2}\boldsymbol{%
\varepsilon _{1}}+c\boldsymbol{n}, \\
\nabla _{1}\boldsymbol{n} &=&-a\boldsymbol{\varepsilon _{1}},\ \ \ \nabla
_{2}\boldsymbol{n}=-c\boldsymbol{\varepsilon _{2}},
\end{eqnarray*}

We denote by $\widetilde{\nabla}_{1}f$ and $\widetilde{\nabla}_{2}f$ the
derivatives of Pfaff of $f$ along the curves $\omega_{32}=0, \omega_{31}=0$
respectively. One finds

\begin{equation}
\widetilde{\nabla}_{1}f =-\frac{1}{a}\nabla_{1}f,\ \ \ \widetilde{\nabla}%
_{2}f =-\frac{1}{c}\nabla_{2}f.  \notag
\end{equation}

It follows that

\begin{equation}  \label{3}
\widetilde{\nabla}_{1}\boldsymbol{x}=-\frac{1}{a}\boldsymbol{\varepsilon_{1}}%
,\ \ \ \widetilde{\nabla}_{2}\boldsymbol{x}=-\frac{1}{c}\boldsymbol{%
\varepsilon_{2}},
\end{equation}
\begin{equation}  \label{4}
\widetilde{\nabla}_{1}\boldsymbol{\varepsilon_{1}}=p_{1}\boldsymbol{%
\varepsilon_{2}}-\boldsymbol{n},\ \ \ \widetilde{\nabla}_{2}\boldsymbol{%
\varepsilon_{1}}=p_{2}\boldsymbol{\varepsilon_{2}},
\end{equation}
\begin{equation}  \label{5}
\widetilde{\nabla}_{1}\boldsymbol{\varepsilon_{2}}=-p_{1}\boldsymbol{%
\varepsilon_{1}},\ \ \ \ \widetilde{\nabla}_{2}\boldsymbol{\varepsilon_{2}}%
=-p_{2}\boldsymbol{\varepsilon_{1}}-\boldsymbol{n},
\end{equation}
\begin{equation}
\widetilde{\nabla}_{1}\boldsymbol{n}=\boldsymbol{\varepsilon_{1}},\ \ \
\widetilde{\nabla}_{2}\boldsymbol{n}=\boldsymbol{\varepsilon_{2}},  \notag
\end{equation}
where $p_{1}=-\frac{1}{a}q_{1},p_{2}=-\frac{1}{c}q_{2}$ are the geodesic
curvatures of the spherical curves $\omega_{32}=0$ and $\omega_{31}=0$
respectively. The Mainardi-Codazzi equations have the following forms

\begin{equation}  \label{6}
\widetilde{\nabla}_{1}\frac{1}{c}=p_{2}\Big(\frac{1}{a}-\frac{1}{c}\Big),\ \
\ \widetilde{\nabla}_{2}\frac{1}{a}=p_{1}\Big(\frac{1}{a}-\frac{1}{c}\Big).
\end{equation}

Let $f$ be a sufficient differentiable function on $S$. Then the second
differential parameter of Beltrami corresponding to the fundamental form $%
III $ of $S$ is defined by 
\begin{equation}  \label{777}
\Delta ^{III}f= -\widetilde{\nabla}_{1}\widetilde{\nabla}_{1}f-\widetilde{%
\nabla}_{2}\widetilde{\nabla}_{2}f-p_{2}\widetilde{\nabla}_{1}f+p_{1}%
\widetilde{\nabla}_{2}f.
\end{equation}

Applying (\ref{777}) to the position vector $\boldsymbol{x}$, gives

\begin{equation}
\Delta ^{III}\boldsymbol{x}= -\widetilde{\nabla}_{1}\widetilde{\nabla}_{1}%
\boldsymbol{x}-\widetilde{\nabla}_{2}\widetilde{\nabla}_{2}\boldsymbol{x}%
-p_{2}\widetilde{\nabla}_{1}\boldsymbol{x}+p_{1}\widetilde{\nabla}_{2}%
\boldsymbol{x}.  \notag
\end{equation}

From (\ref{3}) we obtain

\begin{equation}  \label{8}
\Delta ^{III}\boldsymbol{x}= \widetilde{\nabla}_{1}\Big(\frac{1}{a}%
\boldsymbol{\varepsilon_{1}}\Big)+\widetilde{\nabla}_{2}\Big(\frac{1}{c}%
\boldsymbol{\varepsilon_{2}}\Big)+\frac{1}{a}p_{2}\boldsymbol{\varepsilon_{1}%
}-\frac{1}{c}p_{1}\boldsymbol{\varepsilon_{2}}.
\end{equation}

Using (\ref{4}) and (\ref{5}), equation (\ref{8}) becomes

\begin{equation}  \label{9}
\Delta ^{III}\boldsymbol{x}= \Big(\widetilde{\nabla}_{1}\frac{1}{a}\Big)%
\boldsymbol{\varepsilon_{1}}+\Big(\frac{1}{a} -\frac{1}{c}\Big)p_{2}%
\boldsymbol{\varepsilon_{1}} + \Big(\widetilde{\nabla}_{2}\frac{1}{c}\Big)%
\boldsymbol{\varepsilon_{2}}+\Big(\frac{1}{a} -\frac{1}{c}\Big)p_{1}%
\boldsymbol{\varepsilon_{2}}- \Big(\frac{1}{a}+\frac{1}{c}\Big)\boldsymbol{n}%
.
\end{equation}

Taking into account the Mainardi-Codazzi equations (\ref{6}), so equation (%
\ref{9}) reduces to

\begin{equation*}
\Delta ^{III}\boldsymbol{x}=\Big(\widetilde{\nabla }_{1}\big(\frac{1}{a}+%
\frac{1}{c}\big)\Big)\boldsymbol{\varepsilon _{1}}+\Big(\widetilde{\nabla }%
_{2}\big(\frac{1}{a}+\frac{1}{c}\big)\Big)\boldsymbol{\varepsilon _{2}}-\big(%
\frac{1}{a}+\frac{1}{c}\big)\boldsymbol{n}
\end{equation*}%
or equivalently, (see \cite{R13})

\begin{equation}  \label{11}
\Delta ^{III}\boldsymbol{x}= grad^{III}\Big(\frac{2H}{K}\Big)-\Big(\frac{2H}{%
K}\Big)\boldsymbol{n}.
\end{equation}

\begin{remark}
S. Stamatakis and H. Al-Zoubi proved in \cite{R13} relation (\ref{11}) by
using tensors calculus.
\end{remark}

From (\ref{11}) the following results were proved in \cite{R13}.

\begin{theorem}
\label{T1} A surface $S$ in $E^{3}$ is of 0-type 1 corresponding to the
third fundamental form if and only if $S$ is minimal.
\end{theorem}

\begin{theorem}
\label{T2} A surface $S$ in $E^{3}$ is of $III$-type 1 if and only if $S$ is
part of a sphere.
\end{theorem}

\begin{corollary}
\label{C1.2} The Gauss map of every surface $S$ in $E^{3}$ is of $III$-type
1. The corresponding eigenvalue is $\lambda =2$.
\end{corollary}

Up to now, the only known surfaces of finite $III$-type in $E^{3}$ are parts
of spheres, the minimal surfaces and the parallel of the minimal surfaces
which are of null $III$-type 2. So the following question seems to be
interesting:

\begin{problem}
Other than the surfaces mentioned above, which surfaces in $E^{3}$ are of
finite $III$-type?
\end{problem}

Another generalization of the above problem is to study surfaces in $E^{3}$
with the position vector $\boldsymbol{x}$ satisfying
\begin{equation}
\Delta ^{III}\boldsymbol{x}=A\boldsymbol{x,}  \label{12}
\end{equation}%
where $A\in \mathbb{Re}^{3\times 3}$.

From this point of view, we also pose the following problem

\begin{problem}
Classify all surfaces in $E^{3}$ with the position vector $\boldsymbol{x}$
satisfying relation (\ref{12}).
\end{problem}

Concerning this problem, in \cite{R14} S. Stamatakis and H. Al-Zoubi studied
the class of surfaces of revolution and they proved that: A surface of
revolution $S$ satisfies (\ref{12}) if and only if $S$ is a catenoid or part
of a sphere. Recently, the same authors in \cite{R1} studied the class of
ruled surfaces and the class of quadric surfaces. In particular, they proved
that helicoids and spheres are the only ruled and quadric surfaces
satisfying (\ref{12}) respectively.

This paper provides the first attempt at the study of finite type families
of surfaces in $E^{3}$ corresponding to the third fundamental form. Our main
result is the following

\begin{theorem}
\label{T3} All tubes in $E^{3}$ are of infinite type.
\end{theorem}

Our discussion is local, which means that we show in fact that any open part
of a tube is of infinite Chen type.

\section{Tubes in $E^{3}$}

Let $C: \boldsymbol{\alpha} =\boldsymbol{\alpha}(\emph{t})$, $\mathit{t}%
\epsilon (a,b)$ be a regular unit speed curve of finite length which is
topologically imbedded in $E^{3}$. The total space $N_{\boldsymbol{\alpha}}$
of the normal bundle of $\boldsymbol{\alpha}((a, b))$ in $E^{3}$ is
naturally diffeomorphic to the direct product $(a,b)\times E^{2}$ via the
translation along $\boldsymbol{\alpha}$ with respect to the induced normal
connection. For a sufficiently small $r>0$ the tube of radius $r$ about the
curve $\boldsymbol{\alpha}$ is the set:
\begin{equation}
T_{r}( \boldsymbol{\alpha})=\{exp_{\boldsymbol{\alpha}(t)}\boldsymbol{u}\mid
\boldsymbol{u}\in N_{\boldsymbol{\alpha}} , \ \ \parallel \boldsymbol{u}%
\parallel =r,\ \ t\in(a,b)\}.  \notag
\end{equation}
Assume that ${\mathbf{t}, \mathbf{h}, \mathbf{b}}$ is the Frenet frame and $%
\kappa$ the curvature of the unit speed curve $\boldsymbol{\alpha} =
\boldsymbol{\alpha}(\emph{t})$. For a small real number $r$ satisfies $0 < r
< min\frac{1}{|\kappa|}$, the tube $T_{r}( \boldsymbol{\alpha})$ is a smooth
surface in $E^{3}$, \cite{R12}. Then a parametric representation of the tube
$T_{r}( \boldsymbol{\alpha})$ is given by

\begin{equation}  \label{eq6}
\digamma:\boldsymbol{x}(t,\varphi)= \boldsymbol{\alpha}+ r \cos\varphi%
\mathbf{h}+ r \sin\varphi\mathbf{b}.
\end{equation}

It is easily verified that the first and the second fundamental forms of $%
\digamma$ are given by
\begin{align*}
I &= \big(\delta^{2}+r^{2}\tau^{2}\big)dt^{2} + 2r^{2}\tau dtd\varphi+
r^{2}d\varphi^{2}, \\
II &= \big(-\kappa\delta\cos\varphi+r\tau^{2}\big)dt^{2} + 2r\tau dtd\varphi
+rd\varphi^{2},
\end{align*}
where $\delta: = (1-r\kappa\cos\varphi)$ and $\tau$ is the torsian of the
curve $\boldsymbol{\alpha}$. The Gauss curvature of $\digamma$ is given by

\begin{equation}  \label{eq7}
K =-\frac{\kappa\cos\varphi}{r\delta}.
\end{equation}

Notice that $\kappa\neq0$ since the Gauss curvature vanishes. The Beltrami
operator corresponding to the third fundamental form of $\digamma$ can be
expressed as follows

\begin{eqnarray}  \label{eq8}
\Delta ^{III} &=&\frac{1}{(\kappa \cos \varphi )^{2}}\Bigg[-\frac{\partial
^{2}}{\partial t^{2}}+2\tau \frac{\partial ^{2}}{\partial t\partial \varphi }%
-(\tau ^{2}+\kappa ^{2}\cos^{2}\varphi)\frac{\partial^{2}}{\partial \varphi^{2}}
\\
& &+\frac{\beta }{\kappa \cos \varphi }\frac{\partial }{\partial t}%
+\left(\tau \acute{}+\kappa^{2}\cos\varphi \sin \varphi -\frac{\tau \beta }{%
\kappa \cos \varphi }\right) \frac{\partial }{\partial \varphi }\Bigg] ,
\notag
\end{eqnarray}

where $\beta: =\kappa \acute{}\cos \varphi +\kappa \tau \sin \varphi $ and $%
\acute{}:=\frac{d}{dt}$.

Before we start of proving our main result, we mention and prove the
following special case of tubular surfaces for later use.

\subsection{Anchor rings}

A tube in the Euclidean 3-space is called an anchor ring if the curve $C$ is
a plane circle (or is an open portion of a plane circle). In this case, the
torsian $\tau$ of $\alpha$ vanishes identically and the curvature $\kappa$
of $\alpha$ is a nonzero constant. Then the position vector $\boldsymbol{x}$
of the anchor ring can be expressed as

\begin{equation}  \label{eq9}
\digamma:\boldsymbol{x}(t,\varphi)= \{(a+ r\cos t) \cos\varphi,(a+ r\cos
t)\sin\varphi, r\sin t \},
\end{equation}

\begin{center}
$a > r,\ \ a \epsilon \mathbb{R}.$
\end{center}

The first fundamental form is

\begin{equation*}
I=r^{2}dt^{2}+(a+r\cos t)^{2}d\varphi ^{2},
\end{equation*}%
while the second is

\begin{equation*}
II=rdt^{2}+(a+r\cos t)\cos td\varphi ^{2}.
\end{equation*}

Hence, the Beltrami operator is given by

\begin{equation}  \label{eq10}
\Delta^{III}= -\frac{\partial^{2}}{\partial t^{2}}+\frac{\sin t}{\cos t}%
\frac{\partial}{\partial t} - \frac{1}{\cos^{2}t}\frac{\partial^{2}}{
\partial\varphi^{2}}.
\end{equation}

Let $x_{1}$ be the first coordinate function of $\boldsymbol{x}$. By virtue
of (\ref{eq10}) one can find

\begin{equation}  \label{eq11}
\Delta^{III}x_{1} = \frac{1}{\cos^{2}t} a\cos\varphi + 2r\cos t\cos\varphi.
\end{equation}

Moreover, by a direct computation, we obtain

\begin{equation}  \label{eq12}
(\Delta^{III})^{2}x_{1} = \left(\frac{2}{\cos^{2}t}-\frac{3}{\cos^{4}t}%
\right)a\cos\varphi+ 4r\cos t\cos\varphi,
\end{equation}

\begin{equation}  \label{eq13}
(\Delta^{III})^{3}x_{1} = \left(\frac{4}{\cos^{2}t}-\frac{42}{\cos^{4}t}+%
\frac{45}{\cos^{6}t}\right)a\cos\varphi+ 8r\cos t\cos\varphi.
\end{equation}

It can be seen that $\Delta ^{III}(\cos t\cos \varphi )=2\cos t\cos \varphi $%
, and for each integer $k>0$, it is easy to see that

\begin{equation}  \label{eq14}
\Delta^{III}\frac{\cos\varphi}{\cos^{k}t} = \left(k^{2}-k-\frac{k^{2}-1}{%
\cos^{2}t}\right)\frac{\cos\varphi}{\cos^{k}t}.
\end{equation}

Thus, by induction, one finds

\begin{equation}  \label{eq15}
(\Delta^{III})^{m}x_{1} = \left(\frac{d_{0,m}}{\cos^{2}t}-\frac{d_{1,m}}{%
\cos^{4}t}+-...+\frac{d_{m-1,m}}{\cos^{2m}t}\right)a\cos\varphi+ 2^{m}r\cos
t\cos\varphi,
\end{equation}%
where $d_{j,m}$ are constants, $j = 1, 2,... , m-1$, and

\begin{equation}
d_{0,m}=2^{m-1}, \ \ d_{m-1,m}= (-1)^{m-1}(2m-1)\overset{m}{\underset{j=1}{%
\prod }}(2j-3)^{2}.  \notag
\end{equation}

Notice that $d_{m-1,m} \neq 0$, for each integer $m \geq 1$. Now, if $%
\digamma$ is of finite type, then there exist real numbers, $c_{1}, c_{2},
... , c_{m}$ such that

\begin{equation}  \label{eq16}
(\Delta ^{III})^{m}\boldsymbol{x}+c_{1}(\Delta ^{III})^{m-1}\boldsymbol{x}%
+...+c_{m-1}\Delta ^{III}\boldsymbol{x}+c_{m}\boldsymbol{x}=\mathbf{0}.
\end{equation}

Since $x_{1}=(a+r\cos t)\cos \varphi $ is the first coordinate of $%
\boldsymbol{x}$, (\ref{eq16}), one gets

\begin{equation}  \label{eq17}
(\Delta ^{III})^{m}x_{1}+c_{1}(\Delta ^{III})^{m-1}x_{1}+...+c_{m-1}\Delta
^{III}x_{1}+c_{m}x_{1}= 0.
\end{equation}

From (\ref{eq11}-\ref{eq13}), (\ref{eq15}) and (\ref{eq17}) we obtain that

\begin{eqnarray}
2^{m}r\cos t\cos\varphi + a\cos\varphi \sum_{j=1}^{m}\frac{d_{j-1,m}}{%
\cos^{2j}t}+2^{m-1}c_{1}r\cos t\cos\varphi &&  \notag \\
+ c_{1}a\cos\varphi \sum_{j=1}^{m-1}\frac{d_{j-1,m-1}}{\cos^{2j}t}%
+...+2c_{m-1}r\cos t\cos\varphi &&  \notag \\
+c_{m-1} \frac{a\cos\varphi}{\cos^{2}t} +c_{m}(a+r\cos t)\cos\varphi &=& 0
\notag
\end{eqnarray}%
which can be rewritten as

\begin{equation}  \label{e17}
\frac{d_{m-1,m}}{\cos^{2m}t}+\frac{1}{\cos^{2m-2}t}F(\cos t)=0,
\end{equation}%
where $F(u)$ is a polynomial in $u =\cos t$ of degree $2m-2$.

This is impossible for any $m \geq 1$ since $d_{m-1,m} \neq 0$.
Consequently, we have the following

\begin{corollary}
\label{C1.1} Every anchor ring in the Euclidean 3-space is of infinite type.
\end{corollary}

\section{Proof of the main theorem}

Applying relation (\ref{eq8}) on the position vector $\boldsymbol{x}$ of (%
\ref{eq6}) gives

\begin{equation}
\Delta ^{III}\boldsymbol{x}= \frac{\beta}{\kappa^{3}\cos^{3}\varphi}\mathbf{t%
}+\left(2r\cos\varphi-\frac{1}{\kappa\cos^{2}\varphi}\right)\mathbf{h}%
+2r\sin\varphi\mathbf{b},  \notag
\end{equation}
which can be rewritten as

\begin{equation}  \label{e9}
\Delta ^{III}\boldsymbol{x}= \frac{\beta}{\kappa^{3}\cos^{3}\varphi}\mathbf{t%
}+\frac{1}{ \kappa^{2}\cos^{2}\varphi}\mathbf{P_{1}}(\cos\varphi,
\sin\varphi),
\end{equation}
where $\mathbf{P_{1}}(u, v)$ is a vector valued polynomial in $u, v$ of
degree 3 with functions in $t$ as coefficients. Moreover, by a long
computation, we obtain

\begin{equation}  \label{e10}
( \Delta ^{III})^{2}\boldsymbol{x}= \frac{\beta^{3}}{\kappa^{7}\cos^{7}%
\varphi}\mathbf{t}+\frac{1}{\kappa^{6}\cos^{6}\varphi}\mathbf{P_{2}}%
(\cos\varphi, \sin\varphi),
\end{equation}
where $\mathbf{P_{2}}(u, v)$ is a vector valued polynomial in $u, v$ of
degree 7 with functions in $t$ as coefficients.

We need the following lemma which can be proved directly by using (\ref{eq8}%
).

\begin{lemma}
\label{L1} For any natural numbers m and n we have
\begin{equation}
\left( \Delta^{III}\frac{\beta^{m}}{(\kappa\cos\varphi)^{n}}\right)= -\frac{%
n(n+2)\beta^{m+2}}{(\kappa\cos\varphi)^{n+4}}+\frac{1}{(\kappa\cos%
\varphi)^{n+3}}P(\cos\varphi,\sin\varphi),  \notag
\end{equation}
\end{lemma}

where $P$ is a polynomial in $u, v$ of degree $n$ + 4 with functions in $t$
as coefficients.

Using lemma \ref{L1} and relation (\ref{eq8}) one finds

\begin{equation}  \label{e11}
( \Delta ^{III})^{\lambda}\boldsymbol{x}= d_{\lambda}\frac{\beta^{2\lambda-1}%
}{(\kappa\cos\varphi)^{4\lambda-1}}\mathbf{t}+\frac{1}{(\kappa\cos%
\varphi)^{4\lambda-2}}\mathbf{P_{\lambda}}(\cos\varphi, \sin\varphi),
\end{equation}%
where

\begin{equation}
d_{\lambda}=(-1)^{\lambda-1}\overset{2\lambda-1}{\underset{j=1}{\prod }}%
(2j-1).  \notag
\end{equation}

It can be seen that $d_{\lambda }\neq 0$, for each natural number $\lambda $%
. Moreover, we have

\begin{equation}  \label{e12}
( \Delta ^{III})^{\lambda+1}\boldsymbol{x}= d_{\lambda+1}\frac{%
\beta^{2\lambda+1}}{(\kappa\cos\varphi)^{4\lambda+3}}\mathbf{t}+\frac{1}{%
(\kappa\cos\varphi)^{4\lambda+2}}\mathbf{P_{\lambda+1}}(\cos\varphi,
\sin\varphi).
\end{equation}

Let $\digamma$ be of finite type. Then there exist real numbers, $c_{1},
c_{2}, ... , c_{\lambda}$ such that

\begin{equation}  \label{eq18}
(\Delta ^{III})^{\lambda+1}\boldsymbol{x}+c_{1}(\Delta ^{III})^{\lambda}%
\boldsymbol{x}+...+c_{\lambda}\Delta ^{III}\boldsymbol{x}=\mathbf{0}.
\end{equation}

Using (\ref{e9}-\ref{e12}), one has

\begin{equation}  \label{eq19}
d_{\lambda+1}\frac{\beta^{2\lambda+1}}{\kappa\cos\varphi}\mathbf{t} = Q_{1}%
\mathbf{t}+ Q_{2}\mathbf{h} +Q_{3}\mathbf{b},
\end{equation}%
where $Q_{i}, i =1, 2, 3$, are polynomials in $u, v$ with functions in $t$
as coefficients.

Now, if $\beta \neq 0.$ From (\ref{eq19}) we find

\begin{equation}  \label{eq20}
d_{\lambda+1}\frac{\beta^{2\lambda+1}}{\kappa\cos\varphi} =
Q_{1}(\cos\varphi,\sin\varphi)
\end{equation}

This is impossible, since $Q_{1}$ is polynomial in $\cos\varphi$ and $%
\sin\varphi$. Assume now $\beta = 0$. Then $\kappa \acute{}=0$ and $%
\kappa\tau=0$ so $\kappa = const.\neq 0$ and $\tau = 0$. Therefore the curve
$C$ is a circle, and so $\digamma$ is anchor ring. Hence, $\digamma$ is of
infinite type according to Corollary (\ref{C1.1}). This completes our proof.



\end{document}